\documentclass[12pt]{article}
\usepackage{amssymb,equation,fullpage,epsfig,equation,enumerat}
\usepackage{graphicx}

\usepackage{ulem,color}

\newcommand{\bemat}[1]{\left[\begin{array}{#1}}
\newcommand{\emat}{\end{array}\right]}
\newcommand{\Eq}[1]{(\ref{#1})}
\newcommand{\Z}{{\mathbb Z}}

\newtheorem{theo}{Theorem}
\newtheorem{lemma}[theo]{Lemma}

\renewcommand{\theequation}{\thesection.\arabic{equation}}

\title{Bifurcations in the regularized Ericksen bar model}

\author{M. Grinfeld \\
\normalsize\small Department of Mathematics, \\
\normalsize\small The University of Strathclyde, \\
\normalsize\small Glasgow, G1 1XH, UK. \\
\normalsize\small{\sffamily m.grinfeld@strath.ac.uk}
\and 
G.~J. Lord \\
\normalsize\small Department of Mathematics and Maxwell Institute,\\
\normalsize\small Heriot-Watt University, \\ 
\normalsize\small Edinburgh, EH14 4AS, UK.\\
\normalsize\small {\sffamily g.j.lord@hw.ac.uk}}
\date{\today}

\begin{document}
\maketitle

\begin{abstract}
\noindent We consider the regularized Ericksen model of an elastic bar on an
elastic foundation on an interval with Dirichlet boundary conditions as a
two-parameter bifurcation problem. We explore, using local bifurcation
analysis and continuation methods, the structure of bifurcations from double
zero eigenvalues. Our results provide evidence in support of M\"uller's
conjecture \cite{Muller} concerning the symmetry of local minimizers of the
associated energy functional and describe in detail the structure of the
primary branch connections that occur in this problem. We give a
reformulation of M\"uller's conjecture and suggest two further conjectures
based on the local analysis and numerical observations. We conclude by
analysing a ``loop'' structure that characterizes $(k,3k)$ bifurcations.
\end{abstract}

{\small 
}


\section{Introduction}

In the late eighties J. M. Ball suggested that an interesting and important 
question in material science would be to understand the {\em dynamical}
creation of microstructure \cite{Ball90}.  A model for creating
microstructure would be a dynamical system with a Lyapunov functional that
does not reach its infimum value on, say, the set of $W^{1,2}_0 (0,1)$
functions, the infimum value being achieved instead by a gradient Young
measure.  Thus one might expect to obtain microstructure dynamically, hoping
that as the Lyapunov functional decreases along trajectories, translates in
time will form a minimizing sequence.

One of the candidates for such a process proposed by Ball {\em et al.}
in \cite{BHJPS} is Ericksen's model of an elastic bar on
an elastic foundation \cite{Eric}. This is given by   
\begin{equation}\label{em}
u_{tt} = \left( W'(u_x) + \beta u_{tx} \right)_x - \alpha u
\end{equation}
on the interval $[0,1]$ with the Dirichlet boundary conditions
\begin{equation}\label{bcs}
u(0,t)=u(1,t)=0.
\end{equation}
Here $u := u(x,t)$ is the lateral displacement of the bar, $\beta$ measures the
strength of viscoelastic effects (the term $\beta u_{xxt}$ provides the
dissipation of energy mechanism) and  $\alpha$ measures the strength of bonding
of the bar to the substrate. $W'(u_x)$ is the (non-monotone) stress/strain
relationship; in what follows we specifically take the double well potential
\begin{equation}\label{w}
W(z) = \frac14 ( z^2-1 )^2.
\end{equation}
It is easily checked that 
\begin{equation}\label{fe1}
E_1= \frac12 \int_0^1 \left[ 
u_t^2 + 2 W(u_x) + \frac{\alpha}{2} u^2 \right]\, dx
\end{equation}
is a Lyapunov function for \Eq{em}. 

Friesecke and McLeod \cite{FM} proved that \Eq{em} admits an uncountable
family of steady states that are energetically unstable but locally
asymptotically stable.  They also showed that initial data evolves, roughly,
to a saw-tooth pattern with the same lap number, (i.e minimum number of
non-overlapping intervals where the pattern is monotone) as the initial
data.  In other words, as Friesecke and McLeod put it in the title of their
paper \cite{FM1}, dynamics is a mechanism preventing the formation of finer
and finer microstructure. These results go some way to explain the earlier
numerical results of Swart and Holmes \cite{SH}. 

M\"uller \cite{Muller} considered the regularized version of the Ericksen
model,
\begin{equation}\label{rem}
u_{tt} = \left( W'(u_x) + \beta u_{tx} -\gamma u_{xxx} \right)_x - \alpha u
\end{equation}
on the interval $[0,1]$ with the double Dirichlet boundary conditions
\begin{equation}\label{2bcs}
u(0,t)=u(1,t)=0; \; \; u_{xx} (0,t)= u_{xx}(1,t)=0.
\end{equation}
The main thrust of M\"uller's sophisticated analysis was to describe the
global minimizer of the associated energy functional,
\begin{equation}\label{fe2}
E_2= \frac12 \int_0^1 \left[ 
u_t^2 + \gamma u_{xx}^2 + 2W(u_x) + \frac{\alpha}{2} u^2 \right]\, dx.
\end{equation}
Before we continue, we need to define periodicity more precisely. Consider a
stationary solution $u_0$ of \Eq{rem}. Take its odd extension to $[1,2]$ and
identify the points $x=0$ and $x=2$. If the resulting function is
$D_{2k}$-periodic on the circle for some $k \in \Z$, we say that $u_0$ is
periodic. Then M\"uller's result is that the global minimizer is a
periodic function with a precisely defined dependence of the period on
$\alpha$ and $\gamma$. He also suggested the following conjecture:

{\bf M\"uller's Conjecture \cite{Muller}:} Local minimizers of $E_2$
are periodic.  

Recently, Yip \cite{Yip} has proved this conjecture for solutions of small
energy where $W$ has the form 
\[
W(p) = (|p|-1)^2. 
\]
In this case many calculations of energy of equilibria can be done 
explicitly. This case, with more general boundary conditions,
was also considered in \cite{TZ,VHRT}. 
Nucleation and ripening in the Ericksen problem with the above form of
free energy density is considered from a more thermodynamical point of
view by Huo and I. M\"uller \cite{HuM}.

Finally, in a related paper \cite{VHR}, 
an extension of Ericksen's model to system of two elastic bars
coupled by springs as a model for martensitic phase transitions
is mainly studied numerically.


The dynamics of the regularized Ericksen bar \Eq{rem} was investigated
in \cite{KH}, where global existence of solutions, existence of a compact
attractor and convergence to equilibria was proved. Furthermore, the case of
$\alpha=0$ was investigated in detail and an almost complete
characterization of the structure of the attractor was given in that case.
A. Novick--Cohen has observed that if $\alpha=0$, the set of stationary
solutions of \Eq{rem} is precisely the same as for the Cahn--Hilliard
equation, which was thoroughly investigated in \cite{GN1,GN2}; more work
exploiting this connection between \Eq{rem} and the Cahn-Hilliard equation
is in preparation \cite{GN3}. In particular, the stationary solutions of
\Eq{rem} for the double Dirichlet boundary conditions correspond to the
Cahn-Hilliard equation with mass zero. As a consequence the bifurcation
diagram of the stationary solutions of
\Eq{rem} with $\alpha=0$ contains only supercritical pitchfork bifurcations
from the trivial solutions; only the branches without internal zeros can be
stable.

Other studies of the dynamics of \Eq{rem} include the work of Vainchtein and
co-workers \cite{Vain1, Vain2, VR}, who considered, in particular,
time-dependent Dirichlet boundary conditions (loading/unloading cycles) in
order to study hysteresis effects. 

In this paper we would like to
\begin{itemize}
\item present some evidence towards verifying M\"uller's conjecture
  and reformulate it;
\item explain how the situation for $\alpha=0$ for \Eq{rem} can be reconciled
  with the result of Friesecke and McLeod for \Eq{em} alluded to above.
\end{itemize}

In very recent related work, Healey and Miller \cite{HM},
have considered the two-dimensional version of the problem with hard loading
on the boundary, using methods of global bifurcation theory and numerical
continuation techniques, concentrating on primary bifurcating branches and
characterizing their symmetry.

We use methods of local bifurcation theory. We start by obtaining the
primary and secondary bifurcation points and presenting the results of
numerical continuation using AUTO \cite{Doedel}. In section
\ref{sec:restab} we apply directly the  Lyapunov-Schmidt theory as
detailed in \cite{GS}; this suggests a mechanism for the restabilization of
unstable solutions. The analysis has uncovered an interesting pattern of
primary branch connections which we analyse in section \ref{loop}. To
conclude we present two further conjectures, these are based on the local
analysis backed up with the numerical observations.

\section{Preliminaries}


We start by reviewing the bifurcation structure of the problem.
As shown in \cite{KH}, the eigenvalues $\nu_k$ of the
linearization of \Eq{rem} around  the trivial solution $u=u_t=0$ satisfy
\begin{equation}  \label{eig}
\nu_k = \frac12 \left(
\beta\pi^2k^2 \pm \sqrt{\beta^2 \pi^4 k^4 - 4 (\gamma\pi^4 k^4 - 
\pi^2 k^2 + \alpha)} \right).
\end{equation}
Hence we have the following lemma.

\begin{lemma}
The eigenvalues of the linearization of \Eq{rem} around the trivial
solution $u=u_t=0$ are generically simple, and pass through zero at points
\begin{equation}\label{ga}
\gamma_{k_i} = \frac{\pi^2 k_i^2-\alpha}{\pi^4 k_i^4}, 
\end{equation}
where the integers $k_i$ are ordered by their distance from the number
\[ 
k^*=\frac{\sqrt{2\alpha}}{\pi}. 
\] 
\end{lemma}

For example, if $\alpha=25$, $k_1=2$, $k_2=3$, $k_3=1$, $k_n=n$ for
$n>3$. In other words at $\alpha=25$, based on the eigenfunctions
$\sin(k\pi x)$, the first bifurcating solution branch has one internal 
zero, the second has two internal zeros, the third has no internal
zeros and the $nth$ has $n-1$ internal zeros.

As we will be working in the $(\alpha,\; 1/\gamma)$ plane, it is
convenient to use \Eq{ga} to define
\begin{equation}\label{G}
\Gamma_k = \{ (\alpha, \; 1/\gamma) \, | \, 
\alpha \in [0,\, \pi^2 k^2], \; 1/\gamma = \pi^4 k^4/(\pi^2 k^2 -\alpha)
\}.
\end{equation}

Thus the curves $\Gamma_k$ are the curves on which the linearization has in
its kernel the eigenfunction $\sin (k\pi x)$. Note that $\Gamma_k$ has the
vertical line $\alpha = k^2 \pi^2$ as an asymptote.  

We can also determine the double zero eigenvalue points.
The curves 
$\Gamma_k$ and $\Gamma_l$ will intersect at a point where 
\[
\alpha \equiv \alpha_{k,l} = \frac{\pi^2 k^2 l^2}{k^2+l^2}.
\]
The corresponding values $\gamma_{k,l}$ can be found from 
\[
\gamma_{k,l} = \frac{\pi^2 k^2 - \alpha_{k,l}}{\pi^4 k^4}.
\]
We call these bifurcation points $(k,l)$ bifurcations. Note that one
can concoct any $(k,l)$ bifurcation point, but never a bifurcation point
of multiplicity higher than two.
From these double zero eigenvalue points curves of secondary
bifurcation points emanate and we denote these by $\Gamma_{k\ell}$.

\subsection{Numerical Evidence}
\label{sec:numerics}
\begin{figure}[here]
\begin{center}
{\bf $\alpha=7.5$ \hspace{6.5cm} $\alpha=33$}\\
\includegraphics[width=0.49\textwidth]{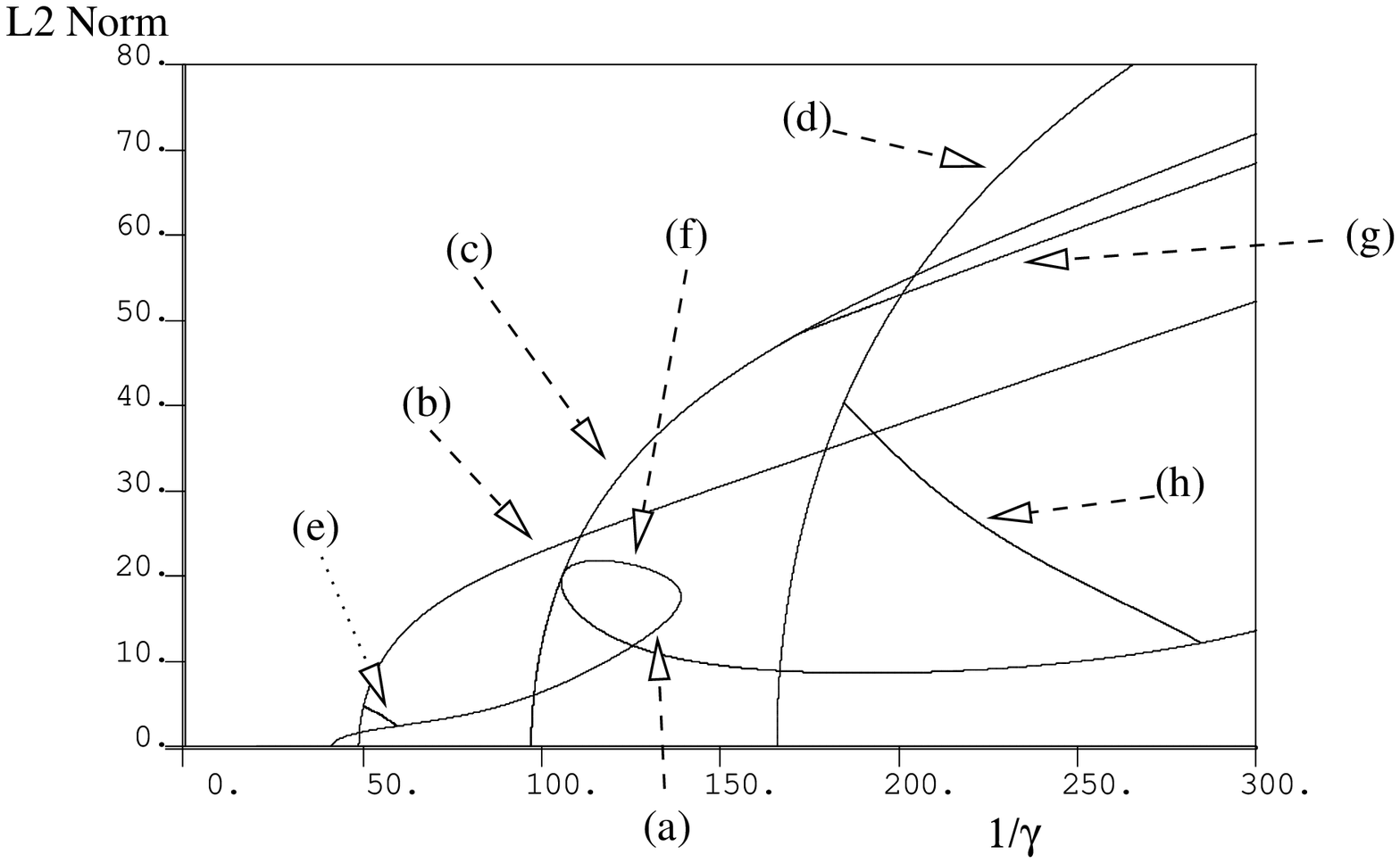}
\includegraphics[width=0.49\textwidth]{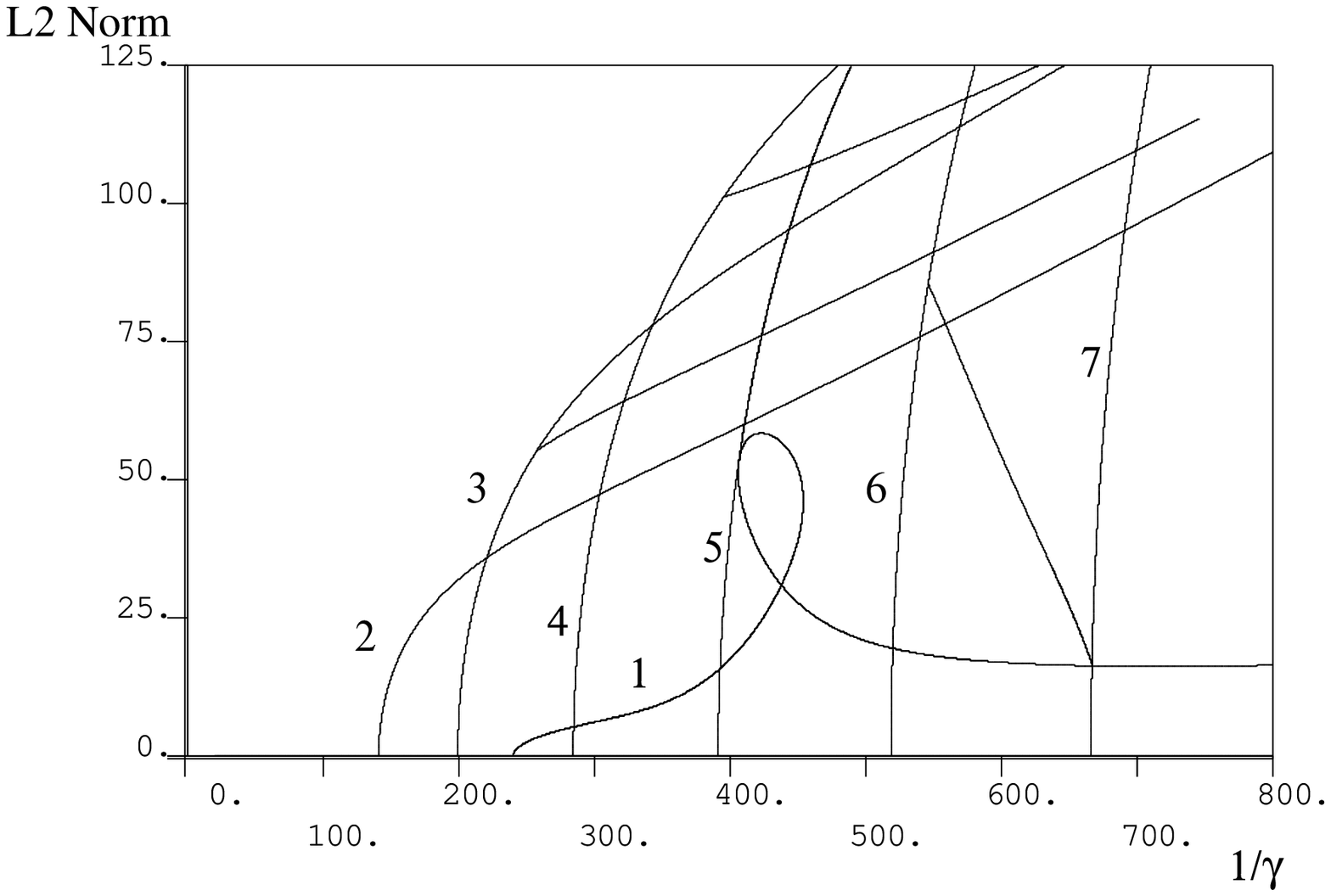}
\caption{Two bifurcation diagrams for fixed values of $\alpha=7.5$
and $\alpha=33$. Figures show the $\ell^2$ norm of 
$(u,u_x,u_{xx},u_{xxx})$ as $\gamma$ varies. In Fig.~\ref{fig:sols} sample
solutions are shown from each of the branches for $\alpha=7.5$; the arrows
give an indication where each the solution is taken from. $(a)$--$(d)$ show 
solution branches with $0, 1, 2$ and $3$ internal zeros. For $\alpha=33$
labels $1$--$7$ show the number of internal zeros.
The ``loop'' for $\alpha=7.5$ connects solutions with $0$ and $2$
internal zeros whereas for $\alpha=33$ it is between solutions with $1$ and $5$ internal
zeros.}
\label{fig:alpha}
\end{center}
\end{figure}

\begin{figure}[here]
\begin{center}
{\bf (a) \hspace{3.2cm} (b) \hspace{3.2cm} (c) \hspace{3.2cm} (d)}\\
\includegraphics[width=0.24\textwidth]{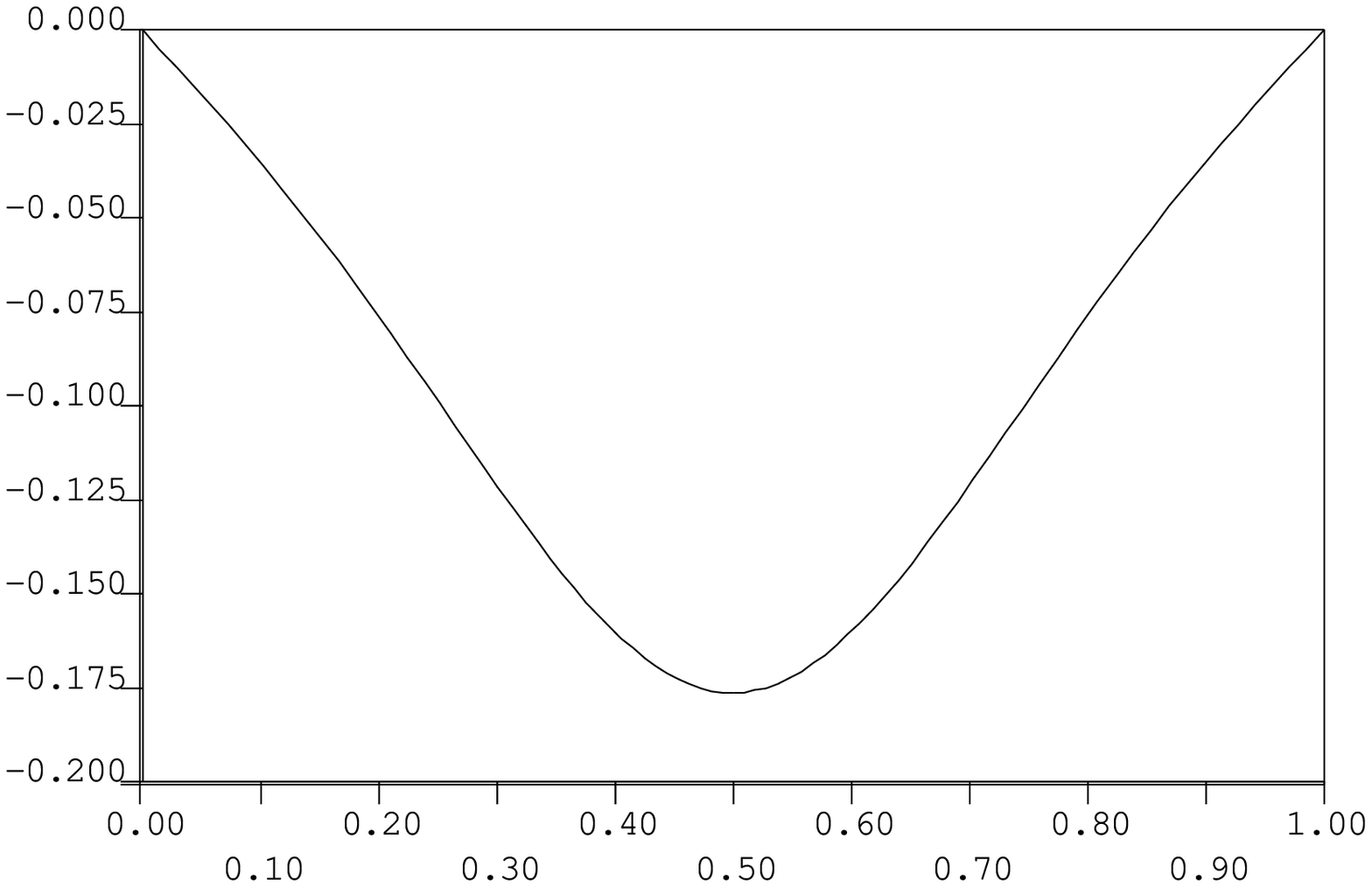}  
\includegraphics[width=0.24\textwidth]{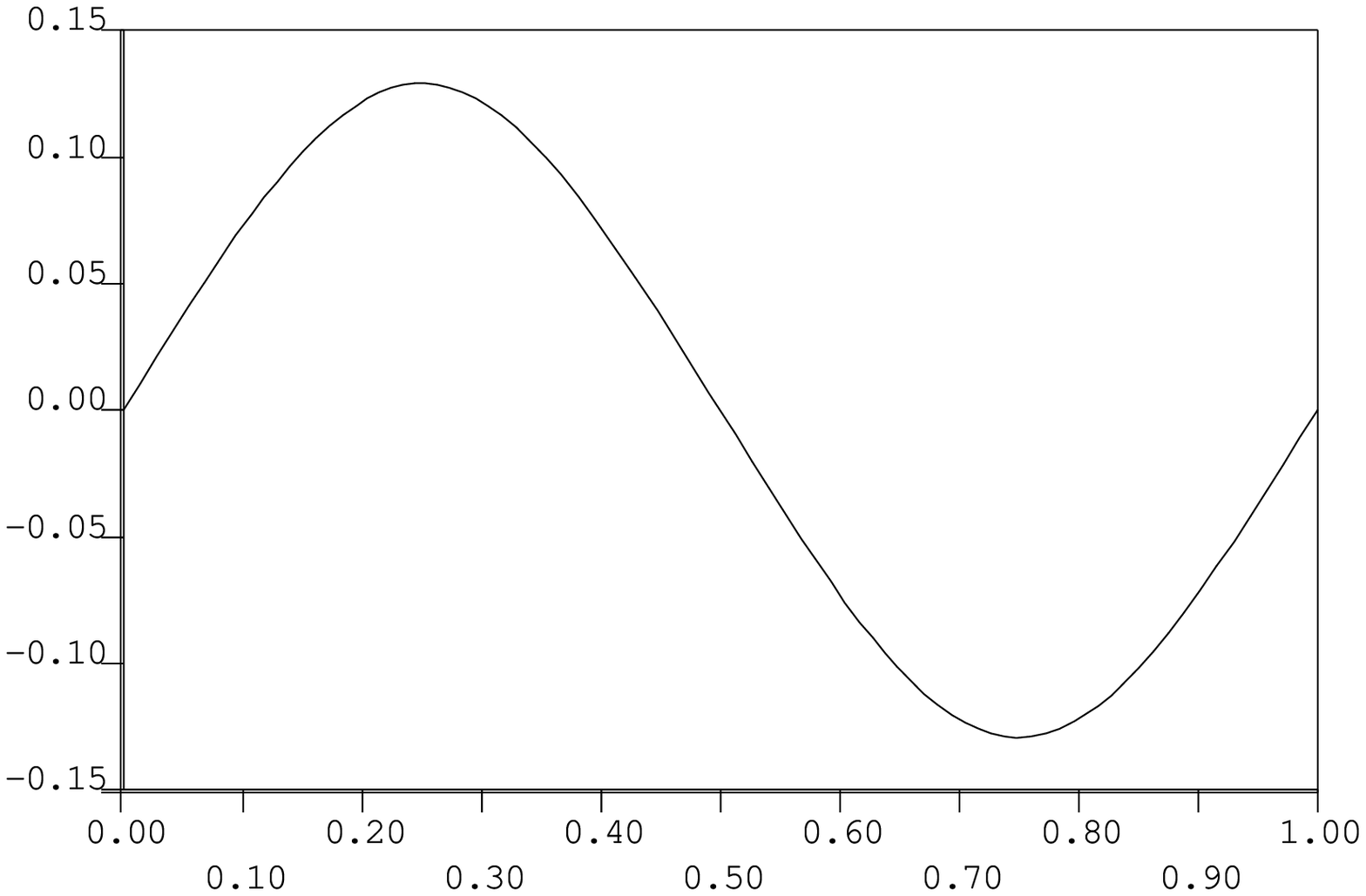}  
\includegraphics[width=0.24\textwidth]{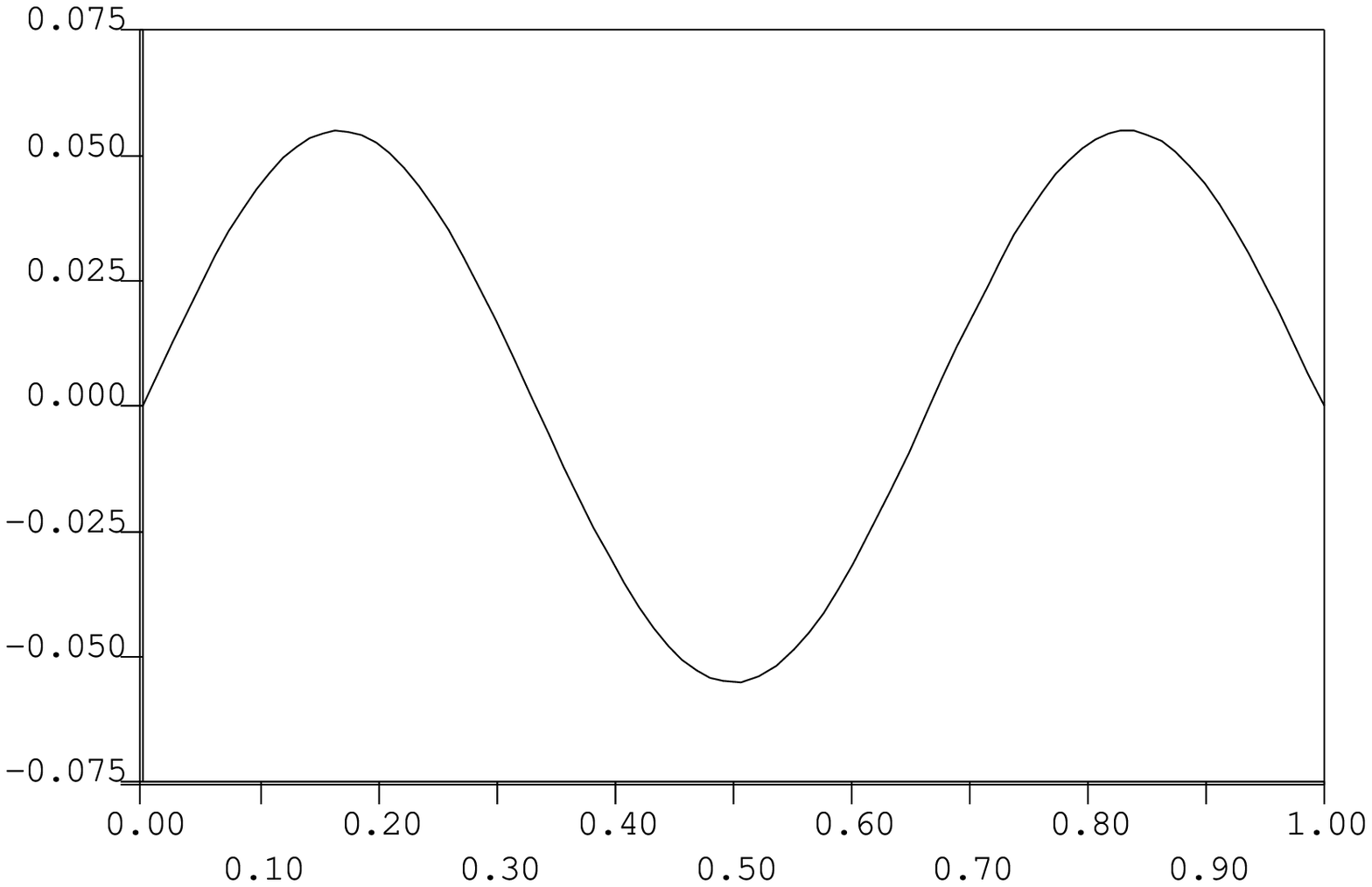}  
\includegraphics[width=0.24\textwidth]{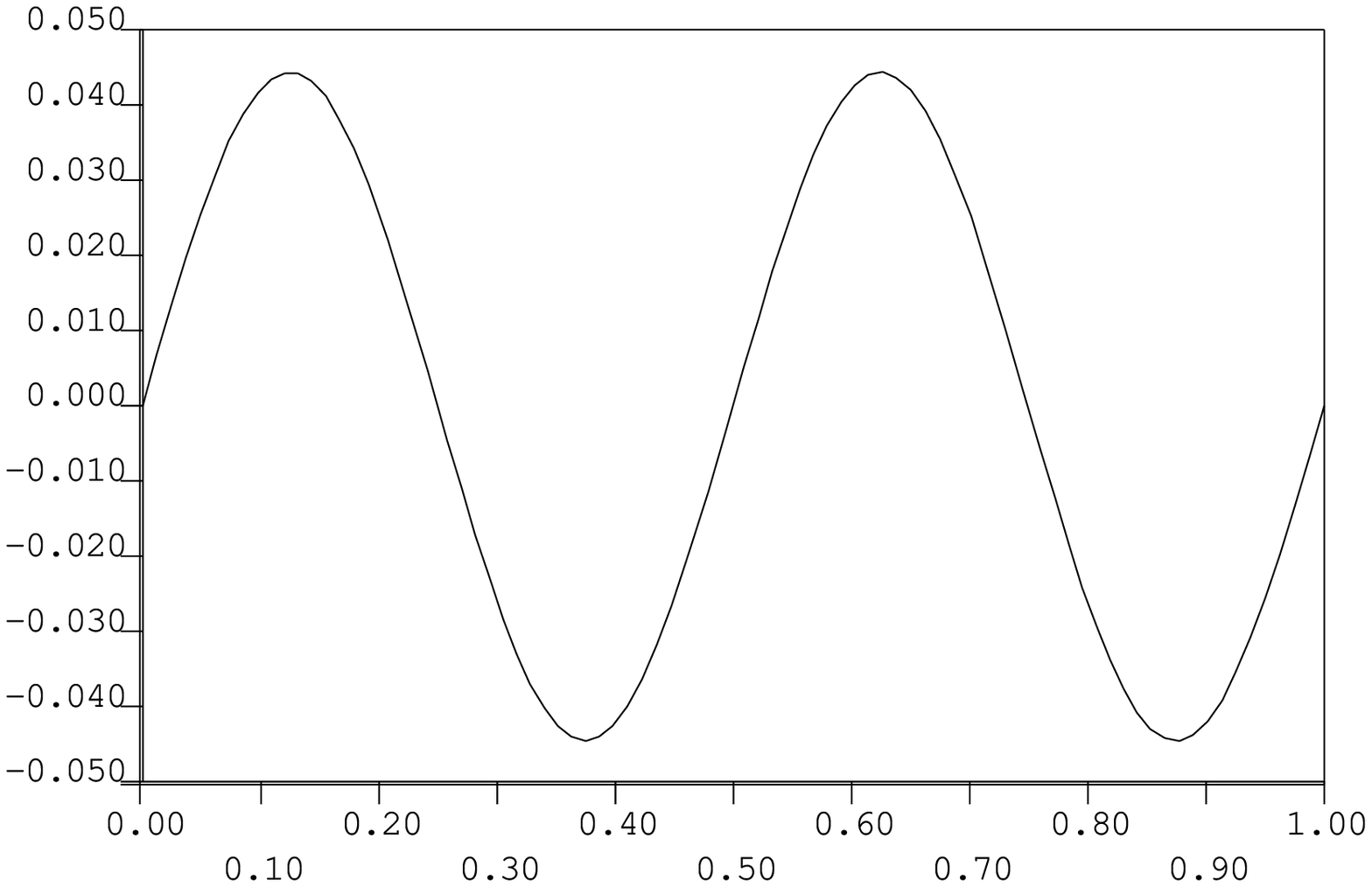}\\
{\bf (e) \hspace{3.2cm} (f) \hspace{3.2cm} (g)  \hspace{3.2cm} (h) }\\
\includegraphics[width=0.24\textwidth]{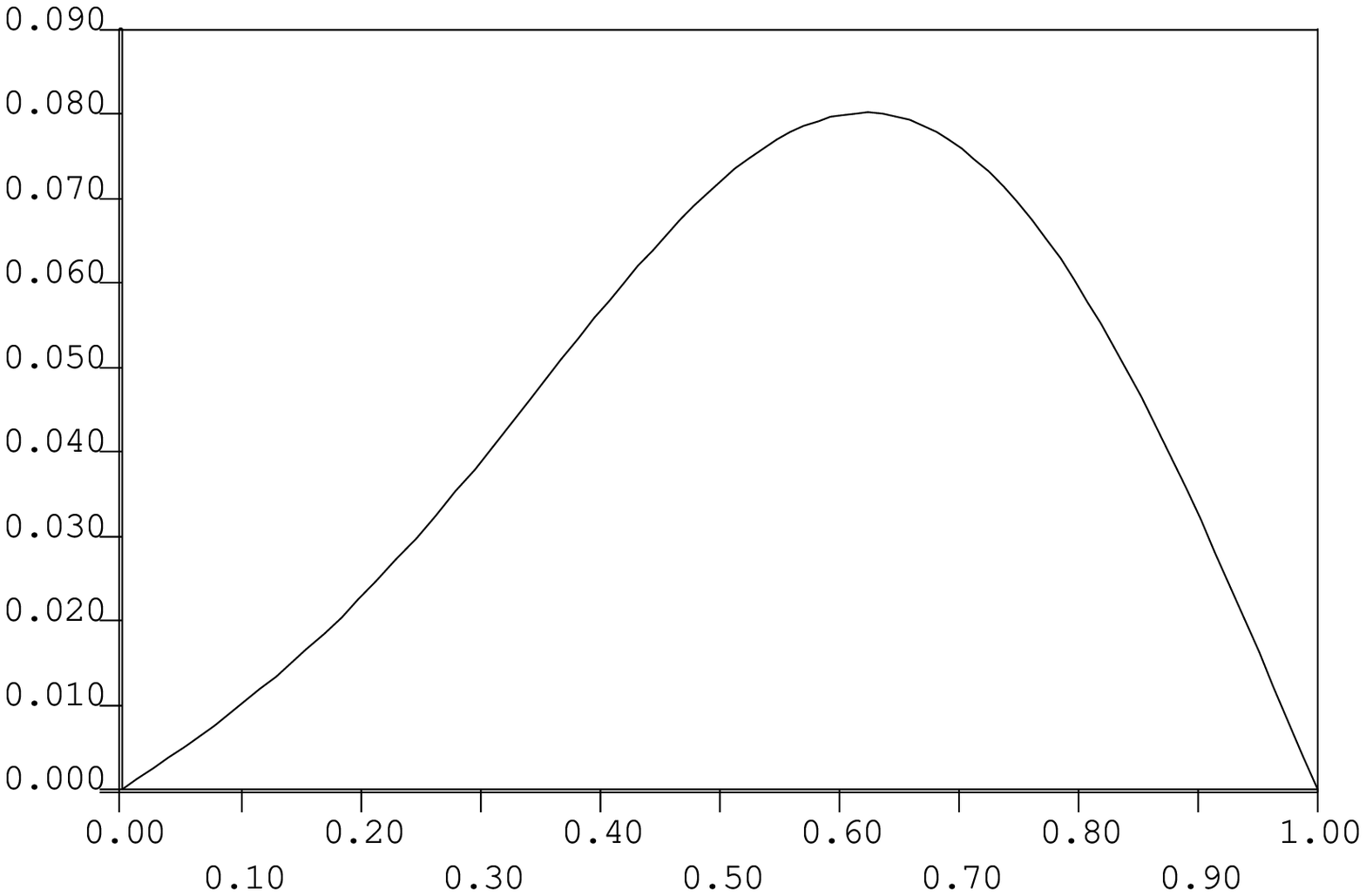}   
\includegraphics[width=0.24\textwidth]{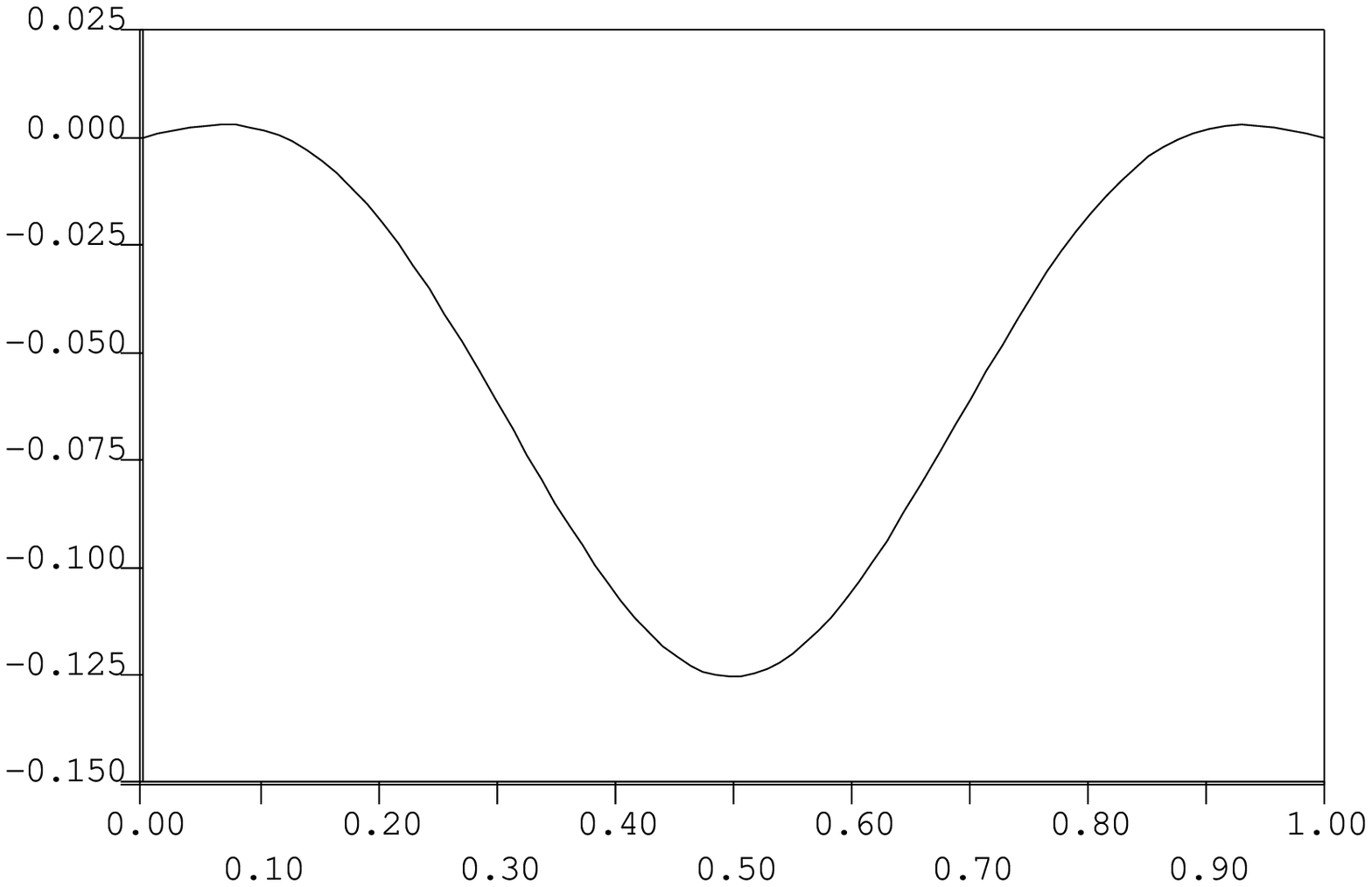}   
\includegraphics[width=0.24\textwidth]{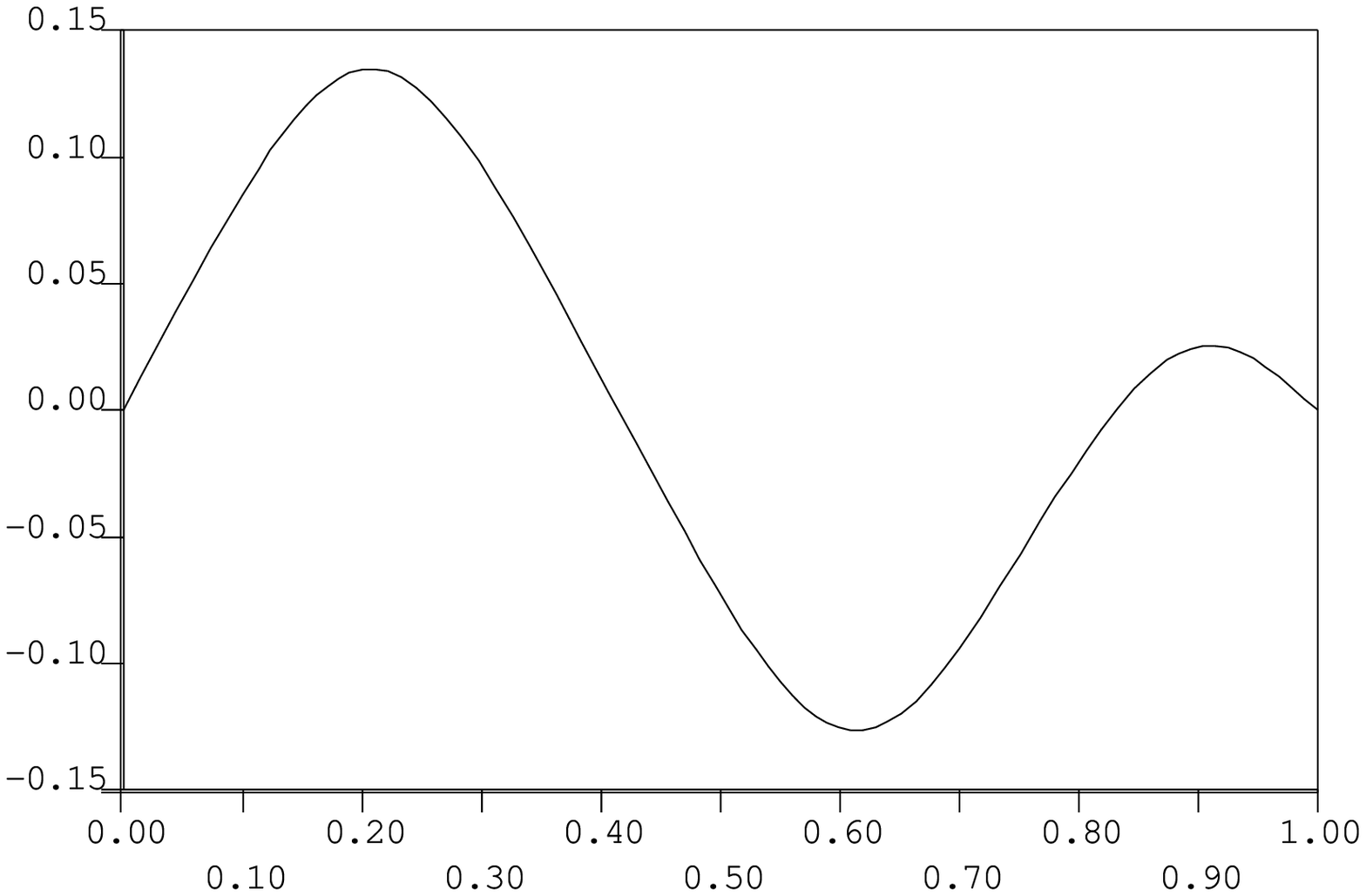} 
\includegraphics[width=0.24\textwidth]{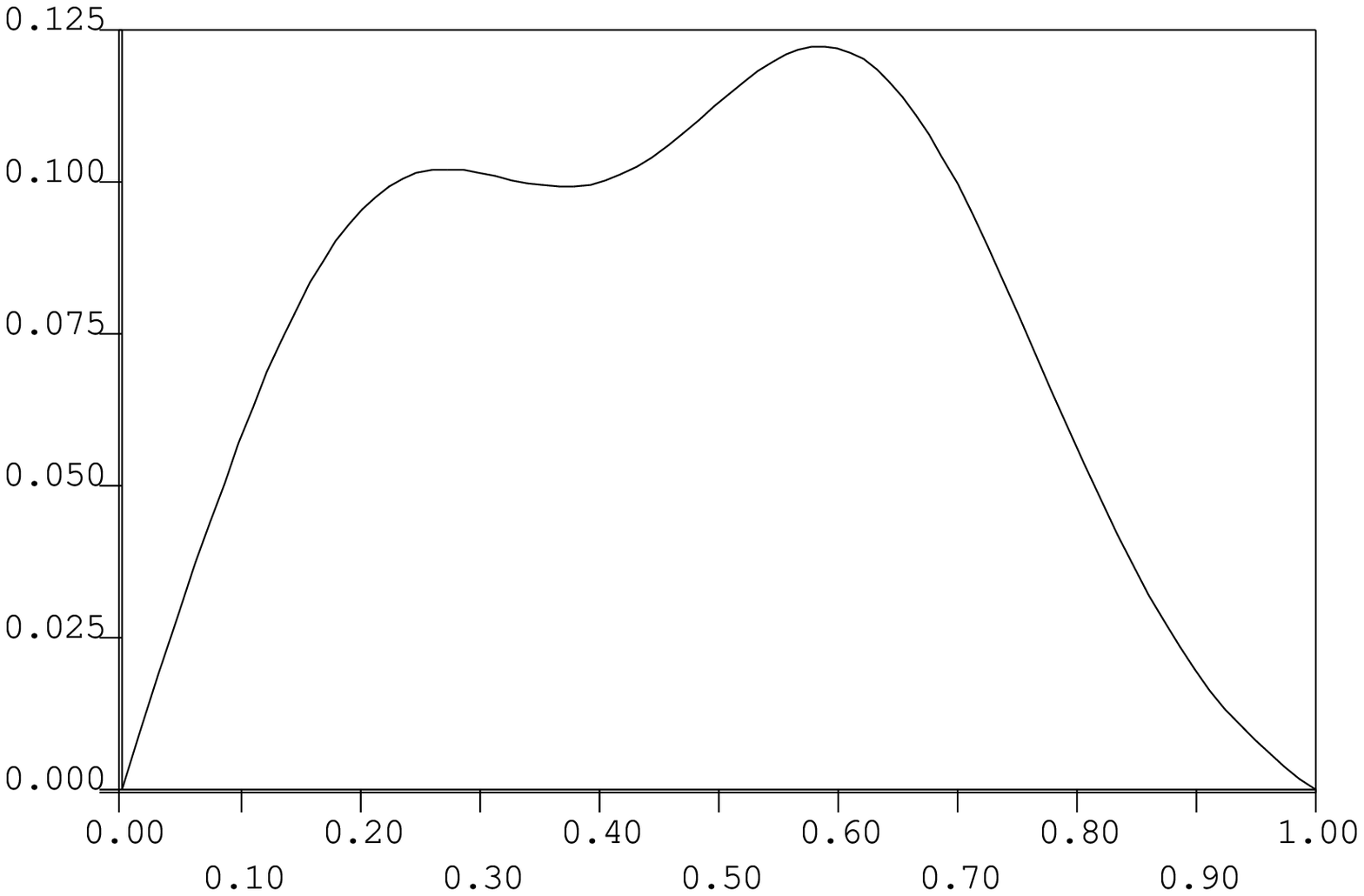}  
\caption{Sample solutions from the bifurcation diagram with $\alpha=7.5$.
{\bf (a)--(d)} shows solutions on the branches bifurcating from the trivial solution with $1/\gamma = 123.87, 111.34,123.78,218.2$.
{\bf (e)--(h)} shows solutions on the secondary bifurcation branches with
$1/\gamma=54.35, 121.40, 278.43, 226.78$.}
\label{fig:sols}
\end{center}
\end{figure}

We use AUTO \cite{Doedel} to investigate numerically the bifurcation
diagram of equilibria of \Eq{rem}, that is we look at
$\Phi(u,\alpha,\gamma) =0$ where 
\begin{equation}\label{phi}
\Phi(u,\alpha,\gamma) =
\gamma u_{xxxx} +\alpha u - (u_x^3-u_x)_x, \;
u(0)=u(1)=0,\qquad  u_{xx}(0)=u_{xx}(1)=0.
\end{equation}
For fixed $\alpha$ we can compute the bifurcation diagram in 
$1/\gamma$ and two examples are shown in Fig.~\ref{fig:alpha} for
$\alpha=7.5$ and $\alpha=33$. The figure shows bifurcations from the
trivial solutions occurring at $\gamma_{k_i}$ and the secondary
bifurcations $(k,\ell)$, plotting the $\ell^2$ norm of 
$(u,u_x,u_{xx},u_{xxx})$ (an approximation of the $H^3$ norm) of the
solution as $\gamma$ varies.
For $\alpha=7.5$, 
$(a)$--$(d)$ are the branches of solutions with $0$--$3$ internal zeros
and sample solutions on these branches are shown in Fig.~\ref{fig:sols}.
Sample solutions from the branches $(e)$--$(h)$ that bifurcate from these
solution branches are also shown in Fig.~\ref{fig:sols}. For
$\alpha=33$ we have labeled the number of internal zeros for the
branches that bifurcate from the trivial solution.

\begin{figure}[here]
\begin{center}
\includegraphics[width=0.9\textwidth]{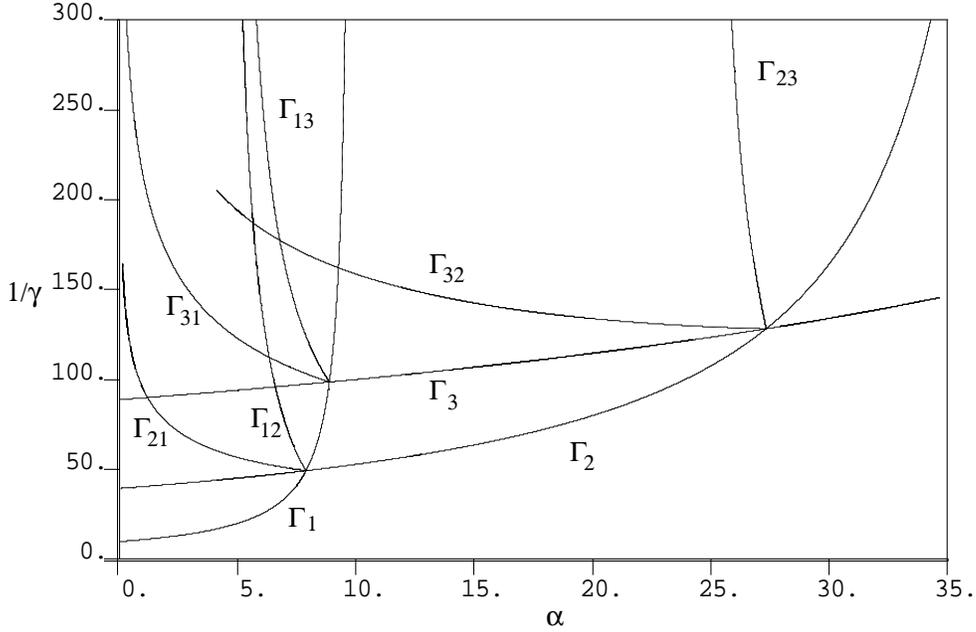}
\caption{Two parameter continuation of the bifurcation points. For
  clarity only the curves $\Gamma_k$ $k=1,2,3$ and the secondary
  curves $\Gamma_{k,\ell}$, $k=1,2,3$, $\ell=1,2,3$ are plotted.}
\label{fig:alphagamma}
\end{center}
\end{figure}

We can exploit the fact that the bifurcations $(k,\ell)$ can be
identified as limit points to perform two parameter continuation. The
result of these computations is shown in the $(\alpha,1/\gamma)$
plane in figure \ref{fig:alphagamma}. 
We clearly see the curve $\Gamma_1$ tending to the correct theoretical
value of the asymptote at $\alpha=\pi^2$. 
From the bifurcation $(1,2)$ the $\Gamma_{21}$ curve
tends to infinite as $\alpha$ approaches zero as does $\Gamma_{31}$;
whereas the curves 
$\Gamma_{12}$, $\Gamma_{13}$, $\Gamma_{23}$ appear to tend to infinite 
for some $\alpha > 0 $.

\section{$(k,k+1)$-bifurcations and a mechanism for restabilization}
\label{sec:restab}

To present a plausible scenario for restabilization of unstable equilibria, 
we are interested in  the structure of stationary solutions of \Eq{phi}
in a neighbourhood of a $(k,k+1)$ bifurcation point, when 
$d\Phi(0,\alpha,\gamma)$ has a double zero eigenvalue with eigenfunctions
$v_k = \sin(k\pi x)$ and $v_{k+1} = \sin((k+1)\pi x)$. 
To examine this we apply the Lyapunov-Schmidt theory as described in
\cite{GS}. 

Set
\[
X = \{ u \in C^4((0,1)), \;| \;  u(0)=u(1)=0, \; u_{xx}(0)=u_{xx}(1)=0.
\},
\]
and let $Y=C^0((0,1))$. We let $L$ denote the linearization : 
$$L \equiv d\Phi(0,\alpha_{k,k+1},\gamma_{k,k+1}).$$

Let us examine the symmetries of (\ref{phi}). Define on $X$ two 
operators, $R_1$ and $R_2$ by $R_1u = -u$ and $R_2 u(x)=u(1-x)$. It is
easily seen that the group $\{ I, R_1, R_2, R_1 R_2\}$ is isomorphic
to $\Z_2 \oplus \Z_2$, that $\Phi$ commutes with this group, i.e.
\[
R_i \Phi(u,\alpha,\gamma)  = \Phi(R_i u,\alpha,\gamma), 
\]
and that if $k$ is even,
\[
R_2 v_k = v_k \hbox{ while } R_2 v_{k+1} = - v_{k+1}. 
\]
Hence the theory of \cite[Chapter X]{GS} is applicable.

In the Lyapunov-Schmidt framework (see \cite[Chapter VII]{GS}),
\[
X = \hbox{sp} \{ v_k,\, v_{k+1}\} \oplus M, 
\]
where $M$ is the orthogonal complement to $\hbox{ker} \,L$. Similarly,
\[
Y= N \oplus \hbox{range}\, L,
\]
where $N$ is the orthogonal complement of the range of $L$.  Since $L$
is self-adjoint, $N = \hbox{ker} L$, so we take $(v_k,\, v_{k+1})$ to be a
basis of both the kernel of $L$ and of $N$.

By the theory of bifurcations with $\Z_2 \oplus \Z_2$ symmetry, 
the bifurcation equation will be of the form
\[{\mathbf g}(x,y,\mu) = {\mathbf 0},\]
where 
\[ {\mathbf g}(x,y,\mu) = 
\bemat{c} g_1(x,y,\mu) \\ g_2(x,y,\mu) \emat =
\bemat{c} Ax^3 + Bxy^2 + a\mu x \\
Cx^2 y +Dy^3 +b \mu y \emat.
\]

From now on we fix $\alpha$ at a point
$\alpha_{k,k+1}$, and do not any longer indicate the dependence on
$\alpha$.  We take our distinguished parameter to be $\lambda =
1/\gamma$ and let it vary through the critical value
$1/\gamma_{k,k+1}$. Clearly,
\[
\frac{\partial^2 g_i}{\partial \lambda \partial x} =
\frac{d\gamma}{d \lambda} \frac{\partial^2 g_i}{\partial \gamma\partial x} =
-\frac{1}{\gamma^2} \frac{\partial^2 g_i}{\partial \gamma \partial x},
\]
which means that the case of positive $g_{i\gamma x}$, say,
corresponds to the case of negative $g_{i\lambda x}$, so we are in
case (A) of \cite[p. 430]{GS}. Then varying $\alpha$ will unfold the degenerate
bifurcation.

By \cite[Appendix 3]{GS} (see also (1.14) in \cite[Chapter VII]{GS}),
and since $d^2 \Phi \equiv 0$ by oddness, we get, for example, for
$g_1$,
\[
\begin{eqalign}
& \frac{\partial^3 g_1}{\partial x^3}~=~ \langle v_k, \;  
d^3 \Phi(0,\lambda)(v_k,\,v_k,\,v_k) \rangle\\
& \frac{\partial^3 g_1}{\partial x \partial y^2}~=~ \langle v_k, \;  
d^3 \Phi(0,\lambda)(v_k,\,v_{k+1},\,v_{k+1}) \rangle\\
& \frac{\partial^2g_1}{\partial x \partial \gamma}~=~ \langle v_k, \;  
d_\gamma \Phi(0,\lambda)(v_k) \rangle,
\end{eqalign}
\]
with similar expressions holding for the partial derivatives of $g_2$.

Note that, for example,
\[
\begin{eqalign}
d^3 \Phi  (v_k, v_k, v_k) = 
\frac{\partial}{\partial t_1} & \frac{\partial}{\partial t_2}
\frac{\partial}{\partial t_3} |_{t_1=t_2=t_3 =0} \, 
\Phi(v_k \sum_{i=1}^3 t_i)\\
=\frac{\partial}{\partial t_1} & \frac{\partial}{\partial t_2}
\frac{\partial}{\partial t_3} |_{t_1=t_2=t_3 =0} \, 
(-(v_k)_x \sum_{i=1}^3 t_i)^3_x = -18 [(v_k)_x]^2 (v_k)_{xx}.
\end{eqalign}
\]
Hence
\[
A = \frac16 \frac{\partial^3 g_1}{\partial x^3} = 3 k^4 \pi^4 \int_0^1 
\sin^4 (k \pi x) \cos^2 (k\pi x)  =
\frac38 k^4 \pi^4.
\]
Similarly,
\[
D = \frac38 (k+1)^4 \pi^4.
\]
For $a$ and $b$ we have
\[
a = \langle v_k,\, (v_k)_{xxxx} \rangle = \frac12 k^4 \pi^4 \hbox{ and }
b = \langle v_{k+1},\, (v_{k+1})_{xxxx} \rangle = \frac12 (k+1)^4 \pi^4.
\]

Now to compute $B$ and $C$: 
\[
\begin{eqalign}
d^3 \Phi (v_k, v_{k+1}, v_{k+1}) = 
\frac{\partial}{\partial t_1} & \frac{\partial}{\partial t_2}
\frac{\partial}{\partial t_3} |_{t_1=t_2=t_3 =0} \, 
\Phi(t_1 v_k + (t_2 + t_3)v_{k+1})\\
= -\frac{\partial}{\partial t_1} & \frac{\partial}{\partial t_2}
\frac{\partial}{\partial t_3} |_{t_1=t_2=t_3 =0} \, 
\left[(v_k)_x t_1 + (t_2+t_3)(v_{k+1})_x\right]^3_x\\
 = - & 6 \left[ (v_k)_{xx} [(v_{k+1})_x]^2 + 2 (v_k)_x(v_{k+1})_{xx}(v_{k+1})_x
 \right].
\end{eqalign}
\]
Hence
\[
\begin{eqalign}
B =  \frac12 &
\frac{\partial^3 g_1}{\partial x \partial y^2} =
3 \int_0^1 k^2 (k+1)^2 \pi^4 \sin^2 (k\pi x) \cos^2 ((k+1)\pi x) \, dx\\
     +\frac12 & \int_0^1 k(k+1)^3 \pi^4 \sin (2k \pi x) \sin (2(k+1)\pi x) dx.
\end{eqalign}
\]
Since the second integral is zero, we have that
\[ 
B = C = 3 k^2 (k+1)^2 \int_0^1 \sin^2 (k\pi x) \cos^2 ((k+1)\pi x) \, dx
= \frac34 \pi^4 (k+1)^2 k^2.
\]
Now we can reduce the bifurcation equation to normal form. First note that
\[
\epsilon_1 = \hbox{sgn} \, (A) = 1, \; \; \epsilon_2 = \hbox{sgn} \, (a)
\hbox{sgn} \, (-1/\gamma^2) = -1,
\] 
\[
\epsilon_3 = \hbox{sgn} \, (D) = 1, \; \; \epsilon_4 = \hbox{sgn} \, (b)
\hbox{sgn} \, (-1/\gamma^2) = -1,
\]
so that indeed we are in case (A) of \cite[p. 430]{GS}. 

By Proposition 2.3 of \cite[p. 424]{GS} the bifurcation diagram is
determined by the modal parameters
\[
m = \left| \frac{b}{Da}\right|B \hbox{ and } n = \left| \frac{a}{Ab}\right|
C.
\]
Now,
\[
m = \frac{(k+1)^4}{k^4} \frac83 \frac{1}{(k+1)^4} \frac34 (k+1)^2 k^2 = 
2 \frac{(k+1)^2}{k^2}
\]
and
\[
n = \frac{k^4}{(k+1)^4} \frac83 \frac{1}{k^4} \frac34 (k+1)^2 k^2 = 
2 \frac{k^2}{(k+1)^2}.
\]
Hence for all $k$, $m \geq 2$, but $n$ can be either smaller or larger than
one. If $m>1$, $n>1$, we are in region (1) of \cite[p. 433]{GS}; if $m>1$, $n<1$, we
are in region (2).

This picture of local bifurcations from a double zero eigenvalue readily
leads to a number of interesting conclusions.

Let us first take the situation where both $m$ and $n$ are larger than one.
The first such bifurcation is the (3,4) one. So assume that $k \geq 3$. 

\begin{figure}[here]
\begin{center}
{\bf (a) \hspace{6.5cm} (b)}\\
\includegraphics[width=0.49\textwidth]{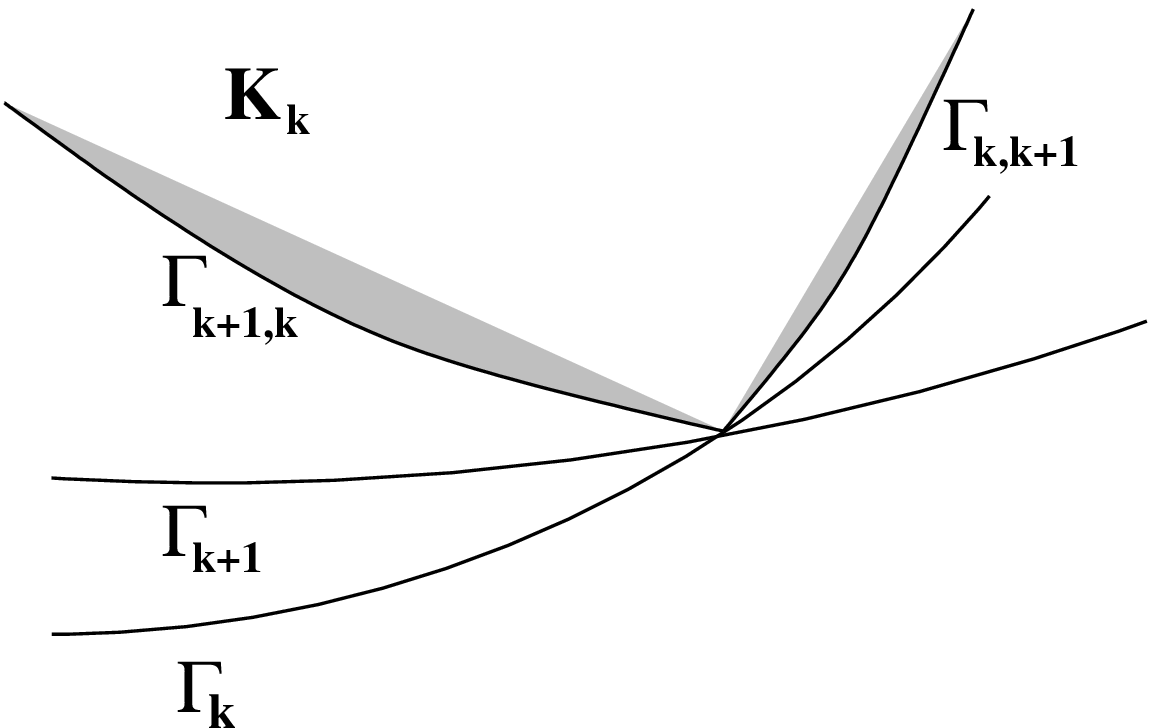}
\includegraphics[width=0.49\textwidth]{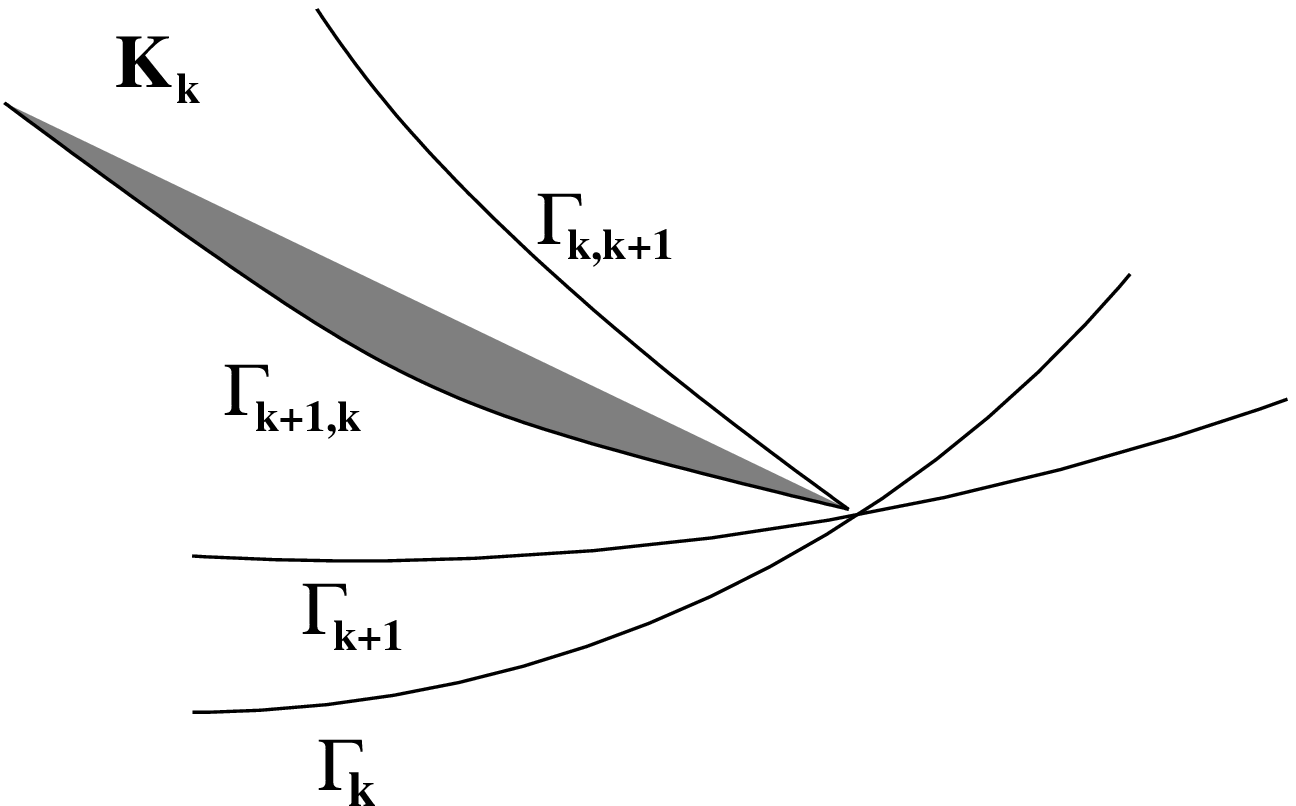}
\caption{Local bifurcation structure as predicted by the
  Lyapunov--Schmidt analysis, in (a) $m,n>1$ and in (b) $m>1, n<1$.
The shaded region indicates the wedge.}
\label{fig:reg}
\end{center}
\end{figure}

Consider the curves $\Gamma_k$ and $\Gamma_{k+1}$ in the $(\alpha,
1/\gamma)$ plane. For $\alpha< \alpha_{k,k+1}$, $\Gamma_k$ lies below
$\Gamma_{k+1}$, and at the bifurcation point $(\alpha_{k,k+1}, 1/\gamma_{k,k+1})$ the
situation reverses. The theory of Golubitsky and Schaeffer tells us that
there are two branches of secondary bifurcations, which we call
$\Gamma_{k+1,k}$ and $\Gamma_{k,k+1}$, from the $(k+1)$-st and the $k$-th
primary curves, respectively, such that close to the double eigenvalue
point $\Gamma_{k,k+1}$ lies to the left of the vertical line
$\alpha=\alpha_{k,k+1}$ and $\Gamma_{k+1,k}$ lies to the right, 
see Fig.~\ref{fig:reg} (a). These two curves of bifurcations form a wedge, 
which we call $K_{k}$. If below the
wedge, to the left of $\alpha_{k,k+1}$ the unstable manifold of solutions on
the $k$-th primary branch has dimension $0$ and that of the $(k+1)$-st
branch has dimension $1$, then to the right of $\alpha_{kl}$ under the
wedge these dimensions switch, while in the interior
of the wedge both the $k$-th and the $(k+1)$-st primary branches are stable. 

In the case of $m>1$, $n<1$, using \cite[Lemma 2.5, p. 426]{GS} we have that
that both $\Gamma_{k+1,k}$ and $\Gamma_{k,k+1}$ lie to the left of the
vertical line $\alpha= \alpha_{k,k+1}$, with $\Gamma_{k,k+1}$ lying above
$\Gamma_{k+1,k}$; as before, this defines a wedge $K_k$ in which both the
$k$-th and the $(k+1)$-st branches are stable.  This case is illustrated 
in Fig.~\ref{fig:reg} (b) and confirms the numerical observation in 
Fig.~\ref{fig:alphagamma}.


\section{$(k,3k)$-bifurcations}\label{loop}

During the process of numerical continuation of branches of solutions, we
discovered an interesting ``loop'' structure connecting each $k$-th and
$3k$-th primary branch (see Fig.~\ref{fig:alpha}). 
We give a simple treatment
based on Fourier mode truncations. We will investigate the $(1,3)$
bifurcation in detail. By a scaling theorem of Aston \cite{Aston}, the
structure of all $(k,3k)$ bifurcations is the same. 

As the required computations are pretty involved, we proceed as follows: (1)
we use $\alpha$ and $\gamma$ as our parameters; this has the effect of
changing hyperbolae into straight lines; (2) we work on the interval
$[0,\pi]$, which allows us to get rid of many multiples of $\pi$; and (3) we
use MAPLE \cite{Maple} throughout, particularly its polynomial
manipulation abilities, such as the computation of resultants and of
Gr\"obner bases, and polynomial factorization. 


Thus, we are looking at
\begin{equation} \label{eq:rescaled}
3u_x^2 u_{xx} -u_{xx} -\gamma u_{xxxx}-\alpha u=0, \qquad  x \in [0,\pi], 
\end{equation}
with double Dirichlet boundary conditions. 

\begin{figure}[here]
\begin{center}
{\bf (a) \hspace{6.5cm} (b)}\\
\includegraphics[width=0.49\textwidth]{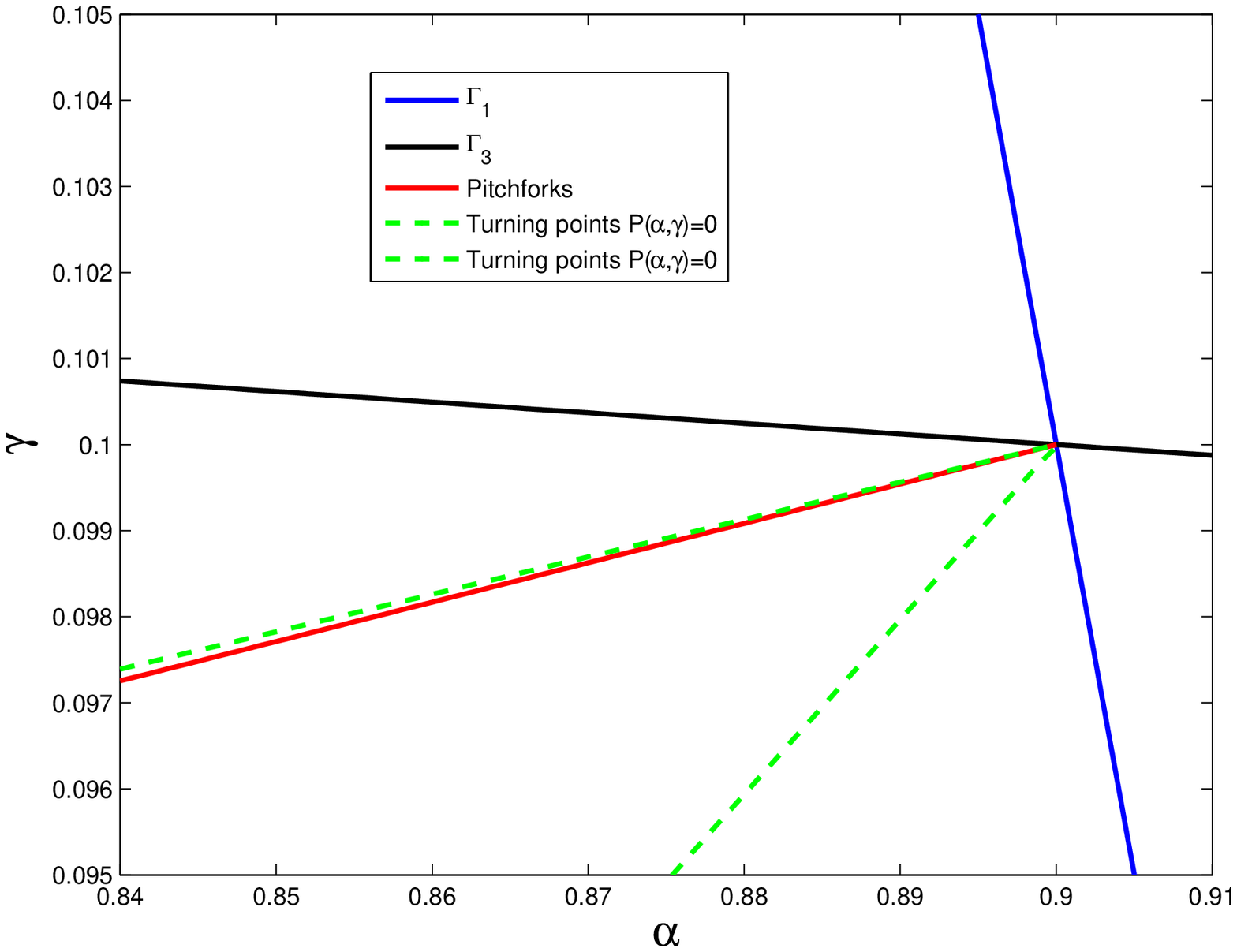}
\includegraphics[width=0.49\textwidth]{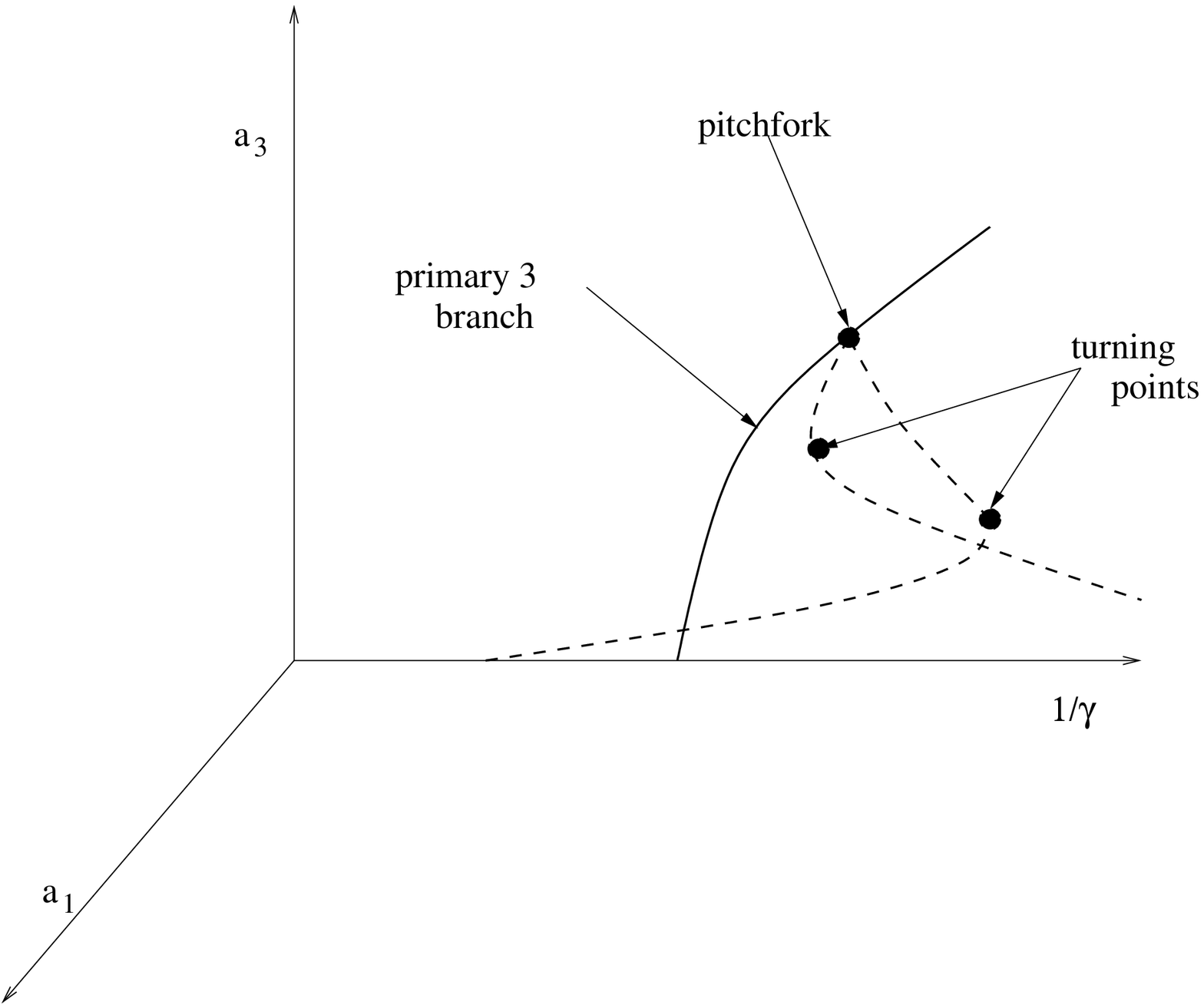}
\caption{{\bf (a)} Plot in the $\alpha$--$\gamma$ plane of the bifurcations 
{\bf (b)} plot illustrating the loop in $a_1$, $a_3$ and $1/\gamma$ space.}
\label{fig:bifk3k}
\end{center}
\end{figure}

Since we are close to the $(1,3)$ bifurcation point, we use the 
ansatz $u= a_1 \sin (x) + a_3 \sin (3x)$. Multiplying  
equation (\ref{eq:rescaled}) by each of $\sin (kx)$, $k=1,\,3$ in turn, 
integrating from 0 to $\pi$ and simplifying,  we get the algebraic
system 
\begin{equation}\label{sys}
\begin{eqalign}
0 & = -\frac12 \alpha a_1- \frac{9}{8} a_1^2 a_3
+\frac12 a_1 -\frac38 a_1^3-1/2\gamma a_1-\frac{27}{4} a_3^{2}a_1\\
0 & =-\frac12\alpha a_3+\frac92 a_3- \frac{243}{8} a_3^3-
\frac{81}{2}\,\gamma a_3-\frac38 a_1^3-\frac{27}{4}a_1^2 a_3
\end{eqalign}.
\end{equation}
Let us examine (\ref{sys}) in more detail.
We need to work out what the curves
$\Gamma_1$ and $\Gamma_3$ look like, and where they intersect. Putting
$a_3$ ($a_1$) to zero in the first (second) of the equations of
(\ref{sys}) gives that  
\[
\Gamma_1 = \{ \gamma= 1-\alpha \}, \qquad \& \qquad \Gamma_3 = \{
\gamma= 1/9 - 1/81 \, \alpha \}. 
\]
These curves are plotted in Fig.~\ref{fig:bifk3k}.
Thus in the $(\alpha, \gamma)$  plane above $\Gamma_1$ there are no non-trivial
solutions of (\ref{sys}). We see that the $(1,3)$ bifurcation point is at
$(\alpha,\gamma)= (9/10,1/10)$ where the two curves intersect.

Now let us see whether (\ref{sys}) will give us a loop, and what this loop
really means. To have a loop we must have a solution with $a_1=0$ and
$a_3\neq 0$  of equations (\ref{sys}).

In other words, we need to solve simultaneously for $a_3$ the equations
\[
-\frac12 \alpha + \frac12 - \frac12\gamma - \frac{27}{4}a_3^2 =0,
\]
and 
\[
-\frac12 \alpha a_3 + \frac92 a_3 -\frac{81}{2}\gamma a_3 
-\frac{243}{8} a_3^3= 0.
\]
These solutions for $a_3$ are the values of $a_3$ on the pure $a_3$ branch
at which a (pitchfork) bifurcation occurs. There are two equations in one
unknown, $a_3$, so for solvability they define a relation between $\alpha$ and
$\gamma$. The form of this relation can be found by taking the
resultant of the above two equations with respect to $a_3$ and is
given by
\begin{equation} \label{pitch}
\gamma = \frac{7}{153} \alpha + \frac{1}{17}.
\end{equation}
This straight line of course intersects the lines $\Gamma_1$ and $\Gamma_3$
at $(9/10,1/10)$ and as the slope is positive, such bifurcation points only
exist for $\alpha<9/10$. For example, it can be checked numerically that if
$\alpha=8/10$, the (secondary) bifurcation from the $a_3$ branch is at
$\gamma = .09542483660$. 

We now turn our attention to the nature of the ``loop''.
For $\alpha < 9/10$ there is no pure $a_1$ branch. The mixed mode branch
bifurcates at $\gamma=1-\alpha$; the pure $a_3$ branch bifurcates off
$\Gamma_3$. It has a secondary bifurcation, a pitchfork, in the $a_1$
direction, on the line (\ref{pitch}). One of these branches then hits the
primary mixed branch and they disappear in a saddle node bifurcation. In
fact, it is possible to work out where the turning points are. They
are given by the following relation connecting $\alpha$ and $\gamma$: 
\begin{equation}\label{hor}
\begin{eqalign}
P_1(\alpha,\gamma) := & 8748\, \gamma+8748\,\alpha+5557\,{\alpha}^{4}-2994732\,\alpha\gamma^{3}\\
+ & 6016437\,\gamma^{4}-
119556\,\alpha \gamma+406782\,{\gamma}^{2}-7938\,{\alpha}^{2}\\
+ & 1054782\,{\alpha}^{2}{\gamma}^{2}-
168492\,{\alpha}^{3}\gamma+243972\,{\alpha}^{2}\gamma\\
- & 1000188\,\alpha{\gamma}^{2}-3156\,{\alpha}^{3}+
288684\,{\gamma}^{3}-2187=0.
\end{eqalign}
\end{equation}
Equation \Eq{hor} is obtained as follows: first divide the first of
the equations of (\ref{sys}) by $a_1$ and then compute the purely lexicographic
Gr\"obner basis of these two equations using the ordering $a_3 > a_1$. This
results in a basis with two elements for the ideal generated by the
two original equations. One of these basis elements 
is found to be a polynomial in $a_1$ only. 
We can then take the resultant of this polynomial in $a_1$ with its
derivative to obtain an enormous polynomial in $\alpha$ and $\gamma$,
eliminating $a_1$.
Some of the roots of this polynomial correspond by construction to
turning points, so it is just a matter of factorizing the resultant
and picking the right term to obtain 
$P_1(\alpha,\gamma)$. If in $P_1(\alpha,\gamma)$ we set $\gamma=1/10-\gamma_1$,
$\alpha=9/10-\alpha_1$, we obtain a homogeneous polynomial of degree 4,
\[
P_2(\alpha_1,\gamma_1) := 5557\alpha_1^4-2994732 \alpha_1 \gamma_1^3+1054782
\alpha_1^2 \gamma_1^2-168492 \alpha_1^3 \gamma_1+6016437 \gamma_1^4.
\]
If we now set $\gamma_1=h \alpha_1$, Descartes' rule of signs tells us that
we will have two real solutions $h_1$ and $h_2$ of
$P_2(\alpha_1,h\alpha_1)/\alpha_1^4 =0$. Using Maple to compute them, we obtain
the approximations for the curves of turning points, 
\[
\gamma \approx 1/10 - 0.203171(9/10-\alpha) \hbox{ and } 
\gamma \approx 1/10 - 0.043484(9/10-\alpha).
\]

In Fig.~\ref{fig:bifk3k} we plot in {\bf (a)} $\Gamma_1$, $\Gamma_2$,
the pitchfork bifurcations given by \Eq{pitch} and the curves of turning
points as determined by \Eq{hor} local to the bifurcation point. We
also illustrate schematically the bifurcations in
$a_1$,$a_3$,$1/\gamma$ space. 

{\bf Remark}: this phenomenon can also be analysed in the framework of
bifurcations with hidden symmetries \cite{AD,Hunt,HM}.

\section{Two new conjectures}

We showed in section \ref{sec:restab} that there are wedges $K_k$ in 
parameter space $(\alpha,1/\gamma)$ where both the $kth$ and $(k+1)st$
branches are stable. Unfortunately though this analysis is only local
in nature. Local to the bifurcation point our numerics of section
\ref{sec:numerics} agrees with the local analysis we have
presented. This lends some confidence to the numerics.
If we combine the local analysis with our observations from the
numerical continuation in Fig.~\ref{fig:alphagamma}, we make the
following conjectures. 

{\bf Conjecture 1.} The wedges $K_k$ are non-empty for all $\gamma <
\gamma_{k,k+1}$. Along the  curves $\Gamma_{k+1,k}$, as $\alpha \rightarrow 0$,
$1/\gamma \rightarrow \infty$, while along $\Gamma_{k,k+1}$, $1/\gamma
\rightarrow \infty$ as $\alpha \rightarrow \widetilde{\alpha_k} \neq 0$.

This automatically means that the curves $\Gamma_{k,k+1}$ and
$\Gamma_{l+1,l}$ intersect if $l>k$.  The upshot of this conjecture is that
the different wedges $K_k$, $K_l$ intersect, and increasingly so as $\gamma
\rightarrow 0$, which for each $N$ creates regions of parameter values in
which there are $N$ different stable equilibria.  This is partially
consistent with the Friesecke and McLeod result, at least for $\alpha <
\widetilde{\alpha_1}$. 
Furthermore we can give a reformulation of {\bf M\"uller's conjecture}
that we stated in the introduction, namely

{\bf M\"uller's conjecture.}
All stable equilibria are created by the above mechanism.

The two statements of this conjecture are equivalent as all the
stabilized branches have $D_{2k}$ symmetry for some $k$.

To conclude we propose another conjecture.

{\bf Conjecture 2.} If $\alpha > k^2 \pi^2$, there are no equilibria with
less than $(k-1)$ internal zeroes. 

This is again based on the local analysis combined with the steady state
bifurcation picture given in section \ref{sec:numerics}, the hope is this
may be easier to prove than the other two. Note that if this conjecture
is true, then for $\alpha > k^2 \pi^2$, an initial condition without
internal zeroes will have to evolve at least $k-1$ interfaces.

\subsection*{Acknowledgements} We are grateful to  A. Novick-Cohen
and G. Berkolaiko for many helpful discussions on this work, and to an
anonymous referee for pointing out much relevant literature, in particular
\cite{HM}.

\end{document}